\documentclass[12pt]{article}
\usepackage{amsmath,amssymb,amsthm,amscd}
\title{Contact resolutions of projectivised nilpotent orbit closures}
\author{Baohua Fu}
\newtheorem{Thm}{Theorem}[section]
\newtheorem{Lem}[Thm]{Lemma}
\newtheorem{Prop}[Thm]{Proposition}
\newtheorem{Cor}[Thm]{Corollary}

\newtheorem{Rque}{Remark}[section]

\def\rit{{\mathbb R}}
\def\cit{{\mathbb C}}

\def\zit{{\mathbb Z}}
\def\pit{{\mathbb P}}

\def\0{{\mathcal O}}
\def\g{{\mathfrak g}}

\def\p{{\mathfrak p}}

\def\n{{\mathfrak n}}

\begin{document}
\maketitle
\begin{abstract}
The projectivised nilpotent orbit  closure $\pit(\overline{\0})$
carries a natural contact structure on its smooth part, which is
induced by a line bundle $L$ on $\pit(\overline{\0})$.  A
resolution $\pi: X \to \pit(\overline{\0})$ is called {\em
contact} if $\pi^* L$ is a contact line bundle on $X$. It turns
out that contact resolutions, crepant resolutions and minimal
models of $\pit(\overline{\0})$ are all the same.
 In this note, we determine when $\pit(\overline{\0})$ admits
a contact resolution, and in the case of existence, we study the
birational geometry among different contact resolutions.
\end{abstract}

\section{Introduction}

Recall that a  nilpotent orbit $\0$ in a semi-simple complex Lie
algebra $\g$ enjoys the following properties:

(i) it is $\cit^*$-invariant, where $\cit^*$ acts on $\g$ by
linear scalars;

(ii) it carries the Kirillov-Kostant-Souriau symplectic 2-form
$\omega$;

(iii) $\lambda^* \omega = \lambda \omega$ for any $\lambda  \in
\cit^*.$

 One deduces from (iii) that this symplectic structure on $\0$
gives a contact structure on the projectivisation $\pit(\0)$,
which is induced by the line bundle
$L:=\0_{\pit(\g)}(1)|_{\pit(\overline{\0})}.$
 When $\g$ is simple, the
variety $\pit(\0) \subset \pit(\g)$ is closed if and only if $\0$
is the minimal nilpotent orbit $\0_{min}$ (see for example Prop.
2.6 \cite{Be}). In this case, $\pit(\0_{min})$ is a Fano contact
manifold. It is generally believed that these are the only
examples of such varieties (\cite{Be}, \cite{Le2}). A positive
answer to this would imply that every compact quaternion-K\"ahler
manifold with positive scalar curvature is homothetic to a Wolf
space (Theorem 3.2 \cite{LS}).

If we take the closure $\overline{\pit(\0)}  =
\pit(\overline{\0})$, then it is in general singular. We say that
a resolution $\pi: X \to \pit(\overline{\0})$ is {\em contact} if
$\pi^*L$ is a contact line bundle on $X$. It follows that $X$ is a
projective contact manifold. Such varieties have drawn much
attention recently(see for example \cite{Pe} and the references
therein).

The first aim of this note is to find all contact resolutions that
$\pit(\overline{\0})$ can have. More precisely we prove that
(Theorem \ref{main}) if the normalization $\pit(\widetilde{\0})$
of $\pit(\overline{\0})$ is not smooth, then the resolution $X$ is
isomorphic to $\pit (T^*(G/P))$ for some parabolic sub-group $P$
in the adjoint group $G$ of $\g$ and $\pi$ is the natural
resolution. The proof relies on the main result in \cite{KPSW} and
that in \cite{Fu}. A classification (Corollary \ref{cla}) of $\0$
such that $\pit(\overline{\0})$ admits a contact resolution can be
derived immediately, with the help of \cite{Be}.

Once we have settled the problem of existence of a  contact
resolution,  we turn to  study the birational geometry among
different contact resolutions in the last section, where (Theorem
\ref{bir}) the chamber structure of the movable cone of a contact
resolution is given, based on the main result in \cite{Na}. This
gives another way to prove the aforesaid result under the
condition that $\overline{\0}$ admits a symplectic resolution,
since minimal models, contact resolutions and crepant resolutions
of $\pit(\overline{\0})$ are the same objects (Proposition
\ref{eq}).
  \vspace{0.3cm}

{\em Acknowledgements:}  The author would like to express his
gratitude to IHES for its hospitality.
 During my visit there, the discussions with
 F. Q. Fang and S. S. Roan are the original impetus for
this work. I am grateful to O. Biquard for discussions on
quaternion-K\"ahler geometry. I want to thank M. Brion for
remarks to a first version of this note, especially
 for the proof  of Proposition \ref{N1}, which ensures the validity
of Theorem \ref{bir} for all simple Lie algebras. I would like to
thank the referees for their critics on a previous version.

\section{Singularities in $\pit(\widetilde{\0})$}

Let $\g$ be a  simple complex Lie algebra and $\0$ a nilpotent
orbit in $\g$. The normalization of the closure $\overline{\0}$
will be denoted by $\widetilde{\0}$. The scalar $\cit^*$-action on
$\overline{\0}$ lifts to $\widetilde{\0}$. There is only one
$\cit^*$-fixed point on $\widetilde{\0}$, say $o$. We denote by
$\pit(\widetilde{\0})$ the geometric  quotient of
 $\widetilde{\0} \setminus \{ o\}$ by the  $\cit^*$-action.
Similarly we denote by $\pit(\overline{\0})$ the geometric
quotient  $\overline{\0} \setminus \{0\}// \cit^*.$ Note that
$\pit(\widetilde{\0})$ is nothing but the normalization of
$\pit(\overline{\0})$.

Recall that a {\em contact structure} on a smooth variety $X$
 is a corank 1
subbundle $F \subset T X$ such that the O'Neil tensor
 $F \times F \to L:=  TX/F$ is everywhere non-degenerate.
In this case, $L$ is called a contact line bundle on $X$ and we
have $K_X \simeq L^{-(n+1)}$, where $n = (dim X -1)/2.$ If we
regard the natural map $T X \to L$ as a section $\theta \in H^0(X,
\Omega_X^1(L))$ (called a {\em contact form}), then the
non-degenerateness is equivalent to the condition that $\theta
\wedge (d \theta)^n$ is nowhere vanishing when considered locally
as an element in $H^0(X, \Omega_X^{2n+1}(L^{n+1})) = H^0(X,
\0_X)$.

For a point $v \in \0$, the tangent space  $T_v \0$ is naturally
isomorphic to  $[v, \g]= Img\ ad_v$.  The map $[v, x] \mapsto
\kappa (v, x) $ defines a one-form $\theta'$ on $\0$, where
$\kappa$ is the Killing form of $\g$. Then $\omega := d \theta'$
is the Kirillov-Kostant-Souriau symplectic form on $\0$. Notice
that $\lambda^* \theta' = \lambda \theta'$ for every $\lambda \in
\cit^*$, so it defines an element $\theta \in H^0(\pit(\0),
\Omega_{\pit(\0)}^1(L)),$ where $L$ is the pull-back of
$\0_{\pit(\g)}(1)$ to $\pit(\0)$. This is in fact a contact form,
i. e. $\theta \wedge (d \theta)^{\wedge n}$
 is everywhere non-zero, where
$n = (dim \0 -2)/2.$  Since the codimension of the complement of
$\pit(\0)$ in $\pit(\widetilde{\0})$ is at least 2, $\theta$
extends to a contact form on the smooth part of
$\pit(\widetilde{\0})$.

\begin{Rque} \upshape
Let $G$ be the adjoint group of $\g$. Then the contact structure
on $\pit(\0)$ is $G$-invariant, which is precisely the contact
structure on $\pit(\0)$ when $\pit(\0)$ is  viewed as a twistor
space of a quaternion-K\"ahler manifold (\cite{Sw}).
\end{Rque}

\begin{Prop}\label{sing}
The projective variety $\pit(\widetilde{\0})$ is projectively normal
with only rational Gorenstein singularities.
\end{Prop}
\begin{proof}
By abusing the notations, we denote also by $L$ the pull-back of
$L$ to the normalization $\pit(\widetilde{\0})$, which is a
 line bundle. Note that the complement of  $\pit(\0)$
in $\pit(\widetilde{\0})$ has codimension at least 2,
 so $K_{\pit(\widetilde{\0})} =
L^{-n-1}$ is locally free, which implies that
$\pit(\widetilde{\0})$ is Gorenstein. Notice that $\widetilde{\0}
\setminus \{o \}$ has rational singularities by results of
Panyushev and Hinich (see \cite{Pa}), so its quotient by the
$\cit^*$ action $\pit(\widetilde{\0})$ has only rational
Gorenstein singularities.
\end{proof}

The following proposition is easily deduced from  Proposition 5.2
in \cite{Be}, which plays an important role to our classification
result (Corollary \ref{cla}).
\begin{Prop} \label{beau}
Let $\g$ be a simple Lie algebra and $\0 \subset \g$ a non-zero
nilpotent orbit. Then $\pit(\widetilde{\0})$ is smooth if and only
if either $\0$ is the minimal nilpotent orbit or $\g$ is of type
$G_2$ and $\0$ is the nilpotent orbit of dimension 8.
\end{Prop}

Singularities in $\pit(\widetilde{\0})$ are examples of the
so-called {\em contact singularities} in \cite{C-F}. Projectivised
nilpotent orbits have already been studied,  for example, in
\cite{Be} (for relation with Fano contact manifolds), \cite{Ko}
(for relation with harmonic maps) and \cite{Sw} (from the twistor
aspect). Their closures have also been studied, for example in
\cite{Po} (for the self-duality), which  give examples of
non-smooth, self-dual projective varieties.

\section{Minimal models}

For a proper morphism between normal varieties $f: X \to W$, we
denote by $N_1(f)$ the vector space (over $\rit$) generated by
reduced irreducible curves contained in fibers of $f$ modulo
numerical equivalence. Let $N^1(f)$ be the group $Pic(X)\otimes
\rit$ modulo numerical equivalence (w. r. t. $N_1(f)$), then we
have a perfect pairing $N_1(f) \times  N^1(f) \to \rit$.

If $f$ is a resolution, then $X$ is called a {\em minimal model}
of $W$ if $K_X$ is $f$-nef.
\begin{Prop}\label{cre}
Let $W$ be a projective normal variety with only canonical
singularities and $f: X \to W$ a resolution. Then $f$ is crepant if
and only if $X$ is a minimal model of $W$.
\end{Prop}
\begin{proof}
If $f$ is crepant, then $K_X = f^* K_W$, which gives
$K_X \cdot [C] = 0$ for every $f$-exceptional curve $C$, so $X$ is a
minimal model of $W$.

Suppose $K_X$ is $f$-nef, then so is $K_X - f^* K_W = \sum_i a_i
E_i$, where $E_i$ are exceptional divisors of $f$. By the
negativity lemma (see Lemma 13-1-4 \cite{Ma}), $a_i \leq 0$ for
all $i$. On the other hand, $W$ has only canonical singularities,
so $a_i \geq 0$, which gives $a_i = 0$ for all $i$, thus  $f$ is
crepant.
\end{proof}

\begin{Cor}
Let $W$ be a projective normal variety with only terminal
singularities and $f: X \to W$ a resolution. Then the following
statements are equivalent:

(i)  $f$ is crepant;

(ii) $X$ is a minimal model of $W$;

(iii) $f$ is small, i.e. $codim (Exc(f)) \geq 2.$
\end{Cor}

By the previous section, there is a contact structure on
$\pit(\0)$, induced by the line bundle $L$ on
$\pit(\overline{\0})$.  A {\em contact resolution} of
$\pit(\widetilde{\0})$ is a resolution $\pi: X \to
\pit(\widetilde{\0})$ such that $\pi^* L$ is a contact line bundle
 on
 $X$.
\begin{Prop} \label{eq}
Let  $\pi: X \to \pit(\widetilde{\0})$ be a resolution, then
the following statements are equivalent:

(i) $\pi$ is crepant;

(ii)  $K_X$ is $\pi$-nef;

(iii)  $\pi$ is a contact resolution.
\end{Prop}
\begin{proof}
The equivalence between (i) and (ii) follows from Prop. \ref{sing}
and Prop. \ref{cre}.  The implication (iii) to (i) is clear from
the definitions. Now suppose that $\pi$ is crepant, then $K_X
\simeq \pi^* (L^{-(n+1)}) \simeq (\pi^*L)^{-(n+1)}.$ Let $\hat{X}$
be the fiber product $X \times_{\pit(\widetilde{\0})}
(\widetilde{\0} \setminus \{o\})$ and $h: \hat{X} \to
\widetilde{\0} \setminus \{o\} $ the natural projection. Then $h$
is a resolution of singularities and
 $h^* \omega$ extends to a 2-form $\tilde{\omega}$ on
$\hat{X}$ since $\widetilde{\0} \setminus \{o\}$ has only
symplectic singularities, where $\omega$ is the symplectic form on
the smooth part of $\widetilde{\0}$. $\hat{X}$ inherits a
$\cit^*$-action from that on $\widetilde{\0}$.
 Contracting $\tilde{\omega}$ with
the vector field generating the $\cit^*$-action, one obtains an
element
$\tilde{\theta} \in H^0(X, \Omega_X \otimes \pi^* L)$. Now it is
clear that $\tilde{\theta}$ gives the contact form on $X$
extending $\theta$.
\end{proof}
\section{Contact resolutions}

Let $f: Z \to \pit(\overline{\0})$  be a resolution and let
$\hat{Z}$ be the fiber product $Z \times_{\pit(\overline{\0})}
W_0$ and $\tilde{f}: \hat{Z} \to W_0$ the natural projection,
where $W_0 = \overline{\0} \setminus \{0\}$. Recall that $L$ is
the restriction of $\0_{\pit(\g)}(1)$ to $\pit(\overline{\0})$.
\begin{Lem} \label{dual}
 $\hat{Z}$ is isomorphic to the
complement of the zero section in the total space of the line
bundle $(f^*L)^*$ and $\tilde{f}$ is a resolution of
singularities.
\end{Lem}
\begin{proof}
This follows from that $W_0$ is naturally isomorphic to the
complement of the zero section in $L^*$ and the fiber product $Z
\times_{\pit(\overline{\0})} L^*$ is isomorphic to $f^*(L^*)
\simeq (f^* L)^*$.
\end{proof}

\begin{Prop}\label{symcont}
The map $f$ is a contact resolution if and only if $\tilde{f}$ is
a symplectic resolution.
\end{Prop}
\begin{proof}
Let $\omega$ be the Kostant-Kirillov-Souriau symplectic form on
$\0$, then $(\tilde{f})^* \omega$ extends to $\tilde{\omega} \in
H^0(\hat{Z}, \Omega_{\hat{Z}}^2)$. $\hat{Z}$ admits a
$\cit^*$-action induced from the one on $W_0$ and for this action,
one has $\lambda^* \tilde{\omega} = \lambda  \tilde{\omega}$ for
all $\lambda \in \cit^*$. By contracting $\tilde{\omega}$ with the
vector field generating the $\cit^*$-action, we obtain a 1-form
$\theta'$ on $\tilde{Z}$ satisfying $\lambda^* \theta' = \lambda
\theta'$, this gives an element $\theta$ in $H^0(Z,
\Omega_Z(f^*L)). $  Then $\theta$ is a contact form if and only if
$\tilde{\omega}$ is a symplectic form.
\end{proof}
From now on, we let $\0$ be a nilpotent orbit such that
$\pit(\widetilde{\0})$ is singular.
\begin{Prop}\label{n}
 Let $\bar{\pi}: X \to
\pit(\overline{\0})$ be a contact resolution and
$\tilde{L}=\bar{\pi}^*(L)$ the contact line bundle on $X$. Then
$(X, \tilde{L})$ is isomorphic to $(\pit(T^*Y),
\0_{\pit(T^*Y)}(1))$ for some smooth projective variety $Y$.
\end{Prop}
\begin{proof}
Note that $K_X \simeq \tilde{L}^{-n-1}$, where $n =(dim \0)/2 -1$.
For a curve $C$ in $X$, we have $K_X \cdot C = -(n+1) L \cdot
\bar{\pi}_*[C]$, thus $K_X$ is not nef. By \cite{KPSW}, $X$ is
either a Fano contact manifold or $(X, \tilde{L})$ is isomorphic
to $(\pit(T^*Y), \0_{\pit(T^*Y)}(1))$ for some smooth projective
variety $Y$.

The map $\bar{\pi}$ factorizes through the normalization, so we
obtain a birational map $\nu: X \to \pit(\widetilde{\0})$.
 By assumption, $\pit(\widetilde{\0})$
is singular. Zariski's main theorem then implies that there exists
a curve $C$ contained in a fiber of $\nu$. Now $K_X \cdot C = 0$,
thus $-K_X$ is not ample, which shows that $X$ is not Fano.
\end{proof}

Let us denote by $\pi_0: \hat{X}  \to  W_0$  the symplectic
resolution provided by Proposition \ref{symcont}. By lemma
\ref{dual}, $\hat{X}$ is isomorphic to $ T^*Y \setminus Y$.
\begin{Lem}
$\pi_0$ extends to a morphism  $\pi: T^*Y \to \overline{\0}$.
\end{Lem}
\begin{proof}
Note that $\pi_0$ lifts to a morphism $\hat{X} \to
\widetilde{W_0}$, where $\widetilde{W_0}$ is the normalization of
$W_0$, which gives a homomorphism $H^0(\widetilde{W_0},
\0_{\widetilde{W_0}}) \to H^0(\hat{X}, \0_{\hat{X}})$.
 Notice that $H^0(\widetilde{W_0}, \0_{\widetilde{W_0}}) = H^0(\widetilde{\0}, \0_{\widetilde{\0}})$ and
$H^0(\hat{X}, \0_{\hat{X}}) = H^0(T^*Y, \0_{T^*Y}).$ On the other
hand, we have a natural morphism $T^*Y \to Spec(H^0(T^*Y,
\0_{T^*Y}))$, which composed with the map $ Spec(H^0(T^*Y,
\0_{T^*Y})) \simeq Spec(H^0(\hat{X}, \0_{\hat{X}})) \to
Spec(H^0(\widetilde{W_0}, \0_{\widetilde{W_0}})) \simeq
Spec(H^0(\widetilde{\0}, \0_{\widetilde{\0}})) = \widetilde{\0}
\to \overline{\0}$ gives  $\pi$.
\end{proof}

Notice that $\pi$ is a symplectic resolution of $\overline{\0}$,
thus the main theorem in \cite{Fu} implies that $\pi$ is
isomorphic to the moment map of the $G$-action on $T^*(G/P)$ for
some parabolic subgroup $P$ in $G$. So we obtain
\begin{Thm} \label{main}
Let $\0$ be a nilpotent orbit in a semi-simple Lie algebra $\g$
such that $\pit(\widetilde{\0})$ is singular.
 Suppose that we have a contact resolution $\pi: Z \to
\pit(\overline{\0})$, then  $Z \simeq \pit(T^*(G/P))$ for some
parabolic subgroup $P$ in the adjoint group $G$ of $\g$ and  the
morphism $\pi$ is the natural one.
\end{Thm}
Now Proposition \ref{beau} implies the following
\begin{Cor}\label{cla}
Suppose  $\g$ is simple.   The projectivised nilpotent orbit
closure $\pit(\overline{\0})$ admits a contact resolution if and
only if either

(i) $\0$ is the minimal nilpotent orbit, or

 (ii) $\g$ is of type $G_2$ and $\0$ is of dimension 8,  or

  (iii) $\overline{\0}$ admits a
symplectic resolution.
\end{Cor}

The classification of nilpotent orbits satisfying case (iii)
 has been done in \cite{Fu} and \cite{Fu2}.
 For example,
every projectivised nilpotent orbit closure in $\mathfrak{sl}_n$
admits a contact resolution, which is given by the
projectivisation of cotangent bundles of some flag varieties.

Recall that the twistor space of a compact quaternion-K\"ahler
manifold is a contact Fano manifold (\cite{Sa}). One may wonder if
a contact resolution of $\pit(\overline{\0})$ could be the twistor
space of a quaterion-K\"ahler manifold. Unfortunately, the answer
to this is in general no, as shown by the following:
\begin{Prop}
Let $G$ be  a simple complex  Lie group with Lie algebra $\g$ and
$P$ a parabolic sub-group of $G$. Then $\pit(T^*(G/P))$ is a
twistor space of a quaternion-K\"ahler manifold if and only if
$G/P \simeq \pit^n$ for some $n$.
\end{Prop}
\begin{proof}
Recall that the image of the moment map $T^*(G/P) \to \g$ is a
nilpotent orbit closure $\overline{\0}$, which gives a generically
finite morphism $\pi: \pit(T^*(G/P)) \to \pit(\overline{\0})$.
There are two cases:

(i) there is a fiber of positive dimension, then as proved in
Proposition \ref{n}, $\pit(T^*(G/P))$ is not Fano.

(ii) every fiber of $\pi$ is zero-dimensional, then $\pi$ is a
finite $G$-equivariant surjective morphism. If $\pit(T^*(G/P))$ is
Fano, then by Proposition 6.3 \cite{Be}, either $\pi = id$ and
$\0= \0_{min}$ or $\pi$ is one of the $G$-covering in the list of
Brylinski-Kostant (see table 6.2 \cite{Be}). In both cases, one
has that $\pit(T^*(G/P))$ is isomorphic to $\pit(\0'_{min})$ for
some minimal nilpotent orbit $\0'_{min} \subset \g',$ which is
possible only if $G/P$ is
 isomorphic to   $\pit^n$ for some $n$.

Now suppose that $\pit(T^*G/P)$ is a twistor space of a
quaternion-K\"ahler manifold $M$.  Then the scalar curvature of
$M$ would be positive, which implies (\cite{Sa}) that
$\pit(T^*G/P)$ is Fano, so $G/P$ is
 isomorphic to   $\pit^n$ for some $n$.
\end{proof}

As pointed out by Prof. A. Swann, this proposition follows also
from \cite{LS}, where it is shown that a contact Fano variety with
$b_2 \geq 2$ is isomorphic to $\pit(T^* \pit^n)$ for some $n$.

\section{Birational geometry}

Let $\g$ be a simple complex Lie algebra and $\0$ a non-zero
nilpotent orbit in $\g$. We now try to understand the birational
geometry between different contact resolutions of
$\pit(\overline{\0})$. We can assume that $\0$ is not the minimal
nilpotent orbit, since $\pit(\overline{\0}_{min})$ is smooth.

Suppose that $\overline{\0}$ admits a symplectic resolution, then
by \cite{Fu}, it is given by the natural  map $\pi: X:= T^*(G/P)
\to \overline{\0}$ for some parabolic
 sub-group $P$ in $G$. Let us denote by $\pi_0$ the restriction
of $\pi$ to $X_0: =T^*(G/P) \setminus (G/P)$, then $\pi_0$ is a
symplectic resolution of $W_0: = \overline{\0} \setminus \{ 0\}. $

I'm indebted to M. Brion for the proof of the following
proposition, which allows us to remove the restriction
that $\g$ is of classical type in an earlier version of this note.

\begin{Prop}\label{N1}
We have $N_1(\pi_0) = N_1(\pi)$ and $N^1(\pi_0) = N^1(\pi)$.
\end{Prop}
\begin{proof}
Consider the natural projections:
 $X_0 \xrightarrow{p_0} G/P \xleftarrow{p} X$, then $Pic(X_0) \otimes \rit$
is identified with $Pic(G/P) \otimes \rit = N^1(G/P)$ via $p_0^*$.
Notice that for
a complete curve $C$ on $X_0$ and a divisor $D \in Pic(G/P)$, we
have $C \cdot p_0^* D = (p_0)_*(C) \cdot D$. Thus we need to show
that images of complete curves in $X_0$ under $(p_0)_*$
 generate $H_2(G/P, \rit) = N_1(G/P)$.

Let $I$ be the set of simple roots of $G$ which are not roots of
the Levi subgroup of $P$, i.e. $I$ is the set of marked roots in
the marked Dynkin diagram of $\p = \text{Lie}(P)$. A basis of
$H_2(G/P, \rit)$ is given by Schubert curves $C_\alpha := P_\alpha
/B$, where $\alpha \in I$ and $P_\alpha$ is the corresponding
minimal parabolic subgroup containing the Borel subgroup $B$.  We
need to lift every $C_\alpha$ to a curve in $X_0$. There are two
cases:

(i) $I$ consists of a single simple root $\alpha$, then $b_2(G/P)
= 1$. Since $\0$ is supposed to be non-minimal, and the
8-dimensional nilpotent orbit closure  in $G_2$  has no symplectic
resolution (Proposition 3.21 \cite{Fu}), by Proposition
\ref{beau},  we can assume that $\widetilde{\0} \setminus \{o\}$
is not smooth. By Zariski's main theorem, there exists a fiber of
$\pi_0$ which has positive dimension. Take an irreducible curve
$C$ containing in  this fiber, then $(p_0)_* C$ is non-zero in
$H_2(G/P, \rit) \simeq \rit$, which generates (over $\rit$)
$N_1(G/P)$.

(ii) $I$ contains at least two simple roots. To lift $C_\alpha$,
we take a simple root $\beta \in I$ different to $\alpha$, then
$\g_\beta$ generates a $G_\alpha$-submodule $M$ of $\g$ contained
in $\n$, where $G_\alpha$ is the simple subgroup of $G$ associated
with the simple root $\alpha$ and $\n$ is the nilradical of $\p$.
Then  in
$T^*(G/P) \simeq G \times^P \n$, there is the closed subvariety
$P_\alpha \times^B M  \simeq G_\alpha \times^{B_\alpha} M $
which is mapped to $G_\alpha M = M$
 with fibers $G_\alpha / B_\alpha \simeq P_\alpha/B,$
where $B_\alpha = B \cap G_\alpha$. Now any fiber of this map
lifts $C_\alpha$.

\end{proof}

Let $\bar{\pi}: \pit(X) \to \pit(\overline{\0})$ be the induced
map, which is
 a contact resolution. The contact structure on $\pit(X)$
is given by the line bundle $\tilde{L} = \0_{\bar{p}}(1)$, where
$\bar{p}: \pit(X) \to G/P$  is the natural map. We have
$Pic(\pit(X)) \simeq Pic(G/P) \oplus \zit [\tilde{L}]$. Notice
that $\tilde{L} = \bar{\pi}^* L$, so for any
$\bar{\pi}$-exceptional curve $C$,  one has $C \cdot \tilde{L} = C
\cdot \bar{\pi}^* L  = 0$, so $\tilde{L}$ is zero in
$N^1(\bar{\pi})$.
 This provides the
identifications $N^1(\bar{\pi}) = N^1(\pi_0)  = N^1(\pi)$ and
$N_1(\bar{\pi}) = N_1(\pi_0)  = N_1(\pi).$

Recall that the cone $NE(\pi) = NE(G/P)$ is generated by Schubert
curves in $G/P$ over $\rit^{\geq 0}$. As shown in the proof of
Proposition \ref{N1}, these Schubert curves are images of curves
in the fibers of $\pi_0$, thus $NE(\pi_0) = NE(\pi)$. Since
$NE(\pi_0) = NE(\bar{\pi})$, we obtain $  NE(\bar{\pi}) =NE(\pi)$.
By  Kleiman's criterion, $\overline{Amp}(\pi_0)
=\overline{Amp}(\pi) =\overline{Amp}(\bar{\pi}). $ By \cite{Na}
Theorem 4.1 (ii), this is a simplicial polyhedral cone.

Let $g: X_0 \to \pit(X)$ and $h: W_0 \to \pit(\overline{\0})$ be
the natural projections, then $\bar{\pi} g = h \pi_0 . $ Take a
$\pi_0$-movable divisor $p_0^* D$, then $(\pi_0)_* p_0^* D = h^*
\bar{\pi}_* \bar{p}^* D \neq 0$, which gives that $ \bar{\pi}_*
\bar{p}^* D \neq 0  $. Notice that $\pi_0^* (\pi_0)_* p_0^* D =
 g^* \bar{\pi}^*  \bar{\pi}_* \bar{p}^* D$ and
$p_0^* D = g^* \bar{p}^* D$, so the cokernel
of $\bar{\pi}^*  \bar{\pi}_* \bar{p}^* D \to \bar{p}^* D $
has support of codimension $\geq 2$. In conclusion
 $\bar{p}^* D$  is $\bar{\pi}$-movable and vice versa.
So we obtain
$\overline{Mov} (\pi_0) =\overline{Mov} (\pi) =
\overline{Mov} (\bar{\pi}).$

To remember the parabolic subgroup $P$, from now on, we will write
$\pi_P$ instead of $\pi$ (similarly for $\pi_0, \bar{\pi}$). For
two parabolic subgroups $Q, Q'$ in $G$, we write $Q \sim Q'$
(called {\em equivalent}) if $T^*(G/Q)$ and $T^*(G/Q')$ give both
symplectic resolutions of a same nilpotent orbit closure. In
\cite{Na}, Namikawa found a way to describe all parabolic
subgroups which are equivalent to a given one. Furthermore the
chamber structure of $\overline{Mov} ({\pi}_P)$ has been described
in {\em loc. cit.} Theorem 4.1. By our precedent discussions
$\overline{Mov} (\pi_0) =\overline{Mov} (\pi) = \overline{Mov}
(\bar{\pi})$, thus we obtain the chamber structure of
$\overline{Mov} (\bar{\pi})$, namely:
\begin{Thm}\label{bir}
Let $\0$ be a non-minimal  nilpotent orbit in a simple complex Lie
algebra $\g$ whose closure admits a symplectic resolution, say
$T^*(G/P)$, where $G$ is the adjoint group of $\g$. Let
$\bar{\pi}_P: \pit(T^*(G/P)) \to \pit(\overline{\0})$ be the
associated contact resolution. Then $\overline{Mov} (\bar{\pi}_P)
= \cup_{Q \sim P}
 \overline{Amp} (\bar{\pi}_Q).$
\end{Thm}

By Mori theory (see for example \cite{Ma}, Theorem 12-2-7), this
implies that every minimal model of $\pit(\overline{\0})$ is of
the form $\pit(T^*(G/Q))$ for some parabolic subgroup $Q \subset
G$ such that $P \sim Q$. Now by Proposition \ref{eq}, this
 gives another proof of Theorem \ref{main} in the case where
$\overline{\0}$ admits a symplectic resolution.

Similarly, as a by-product of our argument, we obtain the
description of the movable cone of a symplectic resolution of
$W_0$, which shows by Mori theory that every symplectic resolution
of $\overline{\0} \setminus \{0\}$ is the restriction of a
Springer map, thus
\begin{Cor}
Let $\g$ be a simple Lie algebra and $\0 \subset \g$ a nilpotent
orbit. Suppose that $\overline{\0}$ admits a symplectic
resolution, then every symplectic resolution of $\overline{\0}
\setminus \{0\}$ can be extended to a symplectic resolution of
$\overline{\0}$.
\end{Cor}

\begin{Rque} \label{o} \upshape
The condition that $\overline{\0}$ admits a symplectic resolution
cannot be removed, due to the following two examples:

(i). if $\g$ is not of type $A$, then $\overline{\0}_{min}$ admits
no symplectic resolution (\cite{Fu}), however $\overline{\0}_{min}
- \{0\}$ is smooth, so trivially it admits a symplectic
resolution;

(ii). if $\g$ is of type $G_2$ and $\0$ is the 8-dimensional
nilpotent orbit, then $W_0: = \overline{\0} - \{0\}$ is not
smooth, and its normalization map $\mu: \widetilde{W}_0 \to W_0$
is a symplectic resolution which does not extends to
$\overline{\0}$, since $\0$ is not a Richardson nilpotent orbit
(Prop. 3.21 \cite{Fu}). Here we used the result in \cite{LSm} and
\cite{Kr} which says that $\widetilde{W}_0$ is in fact the minimal
nilpotent orbit in $\mathfrak{so}_7$, thus it is smooth and
symplectic.

Are these two examples the only exceptions?
\end{Rque}

\quad \\
Laboratoire J. Leray (Math\'ematiques)\\
 Facult\'e  des sciences \\
2, Rue de la Houssini\`ere,  BP 92208 \\
F-44322 Nantes Cedex 03 - France\\
fu@math.univ-nantes.fr
\quad \\
\end{document}